\documentclass{amsart}
\usepackage{amsfonts}
\usepackage{graphicx}
\usepackage{amscd}

\newtheorem{theorem}{Theorem}
\theoremstyle{plain}

\newtheorem{corollary}{Corollary}

\newtheorem{proposition}{Proposition}
\newtheorem{remark}{Remark}

\numberwithin{equation}{section}

\input{tcilatex}

\begin{document}
\title[Ostrowski Type Inequalites]{Some Ostrowski Type Inequalites via
Cauchy's Mean Value Theorem}
\author{S.S. Dragomir}
\address{School of Computer Science and Mathematics\\
Victoria University of Technology\\
PO Box 14428, MCMC 8001\\
VIC, Australia.}
\email{sever@matilda.vu.edu.au}
\urladdr{http://rgmia.vu.edu.au/SSDragomirWeb.html}
\date{February 27, 2003.}
\subjclass{Primary 26D15; Secondary 26D10.}
\keywords{Ostrowski's inequality, Cauchy's mean value theorem.}

\begin{abstract}
Some Ostrowski type inequalities via Cauchy's mean value theorem and
applications for certain particular instances of functions are given.
\end{abstract}

\maketitle

\section{Introduction}

The following result is known in the literature as Ostrowski's inequality 
\cite{AO}.

\begin{theorem}
\label{t1}Let $f:\left[ a,b\right] \rightarrow \mathbb{R}$ be a
differentiable mapping on $\left( a,b\right) $ with the property that $%
\left\vert f^{\prime }\left( t\right) \right\vert \leq M$ for all $t\in
\left( a,b\right) .$ Then 
\begin{equation}
\left\vert f\left( x\right) -\frac{1}{b-a}\int_{a}^{b}f\left( t\right)
dt\right\vert \leq \left[ \frac{1}{4}+\left( \frac{x-\frac{a+b}{2}}{b-a}%
\right) ^{2}\right] \left( b-a\right) M,  \label{1.1}
\end{equation}%
for all $x\in \left[ a,b\right] .$ The constant $\frac{1}{4}$ is best
possible in the sense that it cannot be replaced by a smaller constant.
\end{theorem}

In \cite{SSD1}, the author has proved the following Ostrowski type
inequality.

\begin{theorem}
\label{t2}Let $f:\left[ a,b\right] \rightarrow \mathbb{R}$ be continuous on $%
\left[ a,b\right] $ with $a>0$ and differentiable on $\left( a,b\right) .$
Let $p\in \mathbb{R}\backslash \left\{ 0\right\} $ and assume that 
\begin{equation*}
K_{p}\left( f^{\prime }\right) :=\sup\limits_{u\in \left( a,b\right)
}\left\{ u^{1-p}\left| f^{\prime }\left( u\right) \right| \right\} <\infty .
\end{equation*}
Then we have the inequality 
\begin{multline}
\left| f\left( x\right) -\frac{1}{b-a}\int_{a}^{b}f\left( t\right) dt\right|
\leq \frac{K_{p}\left( f^{\prime }\right) }{\left| p\right| \left(
b-a\right) }  \label{1.2} \\
\times \left\{ 
\begin{array}{ll}
2x^{p}\left( x-A\right) +\left( b-x\right) L_{p}^{p}\left( b,x\right)
-\left( x-a\right) L_{p}^{p}\left( x,a\right) & \text{if }p\in \left(
0,\infty \right) ; \\ 
&  \\ 
\left( x-a\right) L_{p}^{p}\left( x,a\right) -\left( b-x\right)
L_{p}^{p}\left( b,x\right) -2x^{p}\left( x-A\right) & \text{if }p\in \left(
-\infty ,-1\right) \cup \left( -1,0\right) \\ 
&  \\ 
\left( x-a\right) L^{-1}\left( x,a\right) -\left( b-x\right) L^{-1}\left(
b,x\right) -\frac{2}{x}\left( x-A\right) & \text{if }p=-1,%
\end{array}
\right.
\end{multline}
for any $x\in \left( a,b\right) ,$ where for $a\neq b,$%
\begin{equation*}
A=A\left( a,b\right) :=\frac{a+b}{2},\text{ \hspace{0.05in}is the arithmetic
mean,}
\end{equation*}
\begin{equation*}
L_{p}=L_{p}\left( a,b\right) =\left[ \frac{b^{p+1}-a^{p+1}}{\left(
p+1\right) \left( b-a\right) }\right] ^{\frac{1}{p}},\;\text{is the }p-\text{%
logarithmic mean}\;p\in \mathbb{R}\backslash \left\{ -1,0\right\} ,
\end{equation*}
and 
\begin{equation*}
L=L\left( a,b\right) :=\frac{b-a}{\ln b-\ln a}\text{ \hspace{0.05in} is the
logarithmic mean.}
\end{equation*}
\end{theorem}

Another result of this type obtained in the same paper is:

\begin{theorem}
\label{t3}Let $f:\left[ a,b\right] \rightarrow \mathbb{R}$ be continuous on $%
\left[ a,b\right] $ (with $a>0$) and differentiable on $\left( a,b\right) .$
If 
\begin{equation*}
P\left( f^{\prime }\right) :=\sup\limits_{u\in \left( a,b\right) }\left|
uf^{\prime }\left( x\right) \right| <\infty ,
\end{equation*}
then we have the inequality 
\begin{equation}
\left| f\left( x\right) -\frac{1}{b-a}\int_{a}^{b}f\left( t\right) dt\right|
\leq \frac{P\left( f^{\prime }\right) }{b-a}\left[ \ln \left[ \frac{\left[
I\left( x,b\right) \right] ^{b-x}}{\left[ I\left( a,x\right) \right] ^{x-a}}%
\right] +2\left( x-A\right) \ln x\right]  \label{1.3}
\end{equation}
for any $x\in \left( a,b\right) ,$ where for $a\neq b$%
\begin{equation*}
I=I\left( a,b\right) :=\frac{1}{e}\left( \frac{b^{b}}{a^{a}}\right) ^{\frac{1%
}{b-a}},\text{ \hspace{0.05in} is the identric mean.}
\end{equation*}
\end{theorem}

If some local information around the point $x\in \left( a,b\right) $ is
available, then we may state the following result as well \cite{SSD1}.

\begin{theorem}
\label{t4}Let $f:\left[ a,b\right] \rightarrow \mathbb{R}$ be continuous on $%
\left[ a,b\right] $ and differentiable on $\left( a,b\right) .$ Let $p\in
\left( 0,\infty \right) $ and assume, for a given $x\in \left( a,b\right) ,$
we have that 
\begin{equation*}
M_{p}\left( x\right) :=\sup\limits_{u\in \left( a,b\right) }\left\{
\left\vert x-u\right\vert ^{1-p}\left\vert f^{\prime }\left( u\right)
\right\vert \right\} <\infty .
\end{equation*}%
Then we have the inequality 
\begin{multline}
\left\vert f\left( x\right) -\frac{1}{b-a}\int_{a}^{b}f\left( t\right)
dt\right\vert  \label{1.4} \\
\leq \frac{1}{p\left( p+1\right) \left( b-a\right) }\left[ \left( x-a\right)
^{p+1}+\left( b-x\right) ^{p+1}\right] M_{p}\left( x\right) .
\end{multline}
\end{theorem}

For recent results in connection to Ostrowski's inequality see the papers 
\cite{GA1},\cite{GA2} and the monograph \cite{SSD2}.

The main aim of this paper is to point out some generalizations of the
results incorporated in Theorems \ref{t2}-\ref{t4} by the use of Cauchy mean
value theorem. Applications for other particular instances of functions are
given as well.

\section{The Results}

We may state the following theorem.

\begin{theorem}
\label{t2.1}Let $f,g:\left[ a,b\right] \rightarrow \mathbb{R}$ be continuous
on $\left[ a,b\right] $ and differentiable on $\left( a,b\right) .$ If $%
g^{\prime }\left( t\right) \neq 0$ for each $t\in \left( a,b\right) $ and 
\begin{equation}
\left\Vert \frac{f^{\prime }}{g^{\prime }}\right\Vert _{\infty
}:=\sup\limits_{t\in \left( a,b\right) }\left\vert \frac{f^{\prime }\left(
t\right) }{g^{\prime }\left( t\right) }\right\vert <\infty ,  \label{2.1}
\end{equation}%
then for any $x\in \left[ a,b\right] $ one has the inequality 
\begin{multline}
\left\vert f\left( x\right) -\frac{1}{b-a}\int_{a}^{b}f\left( t\right)
dt\right\vert  \label{2.2} \\
\leq \left\vert 2\left( \frac{x-\frac{a+b}{2}}{b-a}\right) g\left( x\right) +%
\frac{\int_{x}^{b}g\left( t\right) dt-\int_{a}^{x}g\left( t\right) dt}{b-a}%
\right\vert \cdot \left\Vert \frac{f^{\prime }}{g^{\prime }}\right\Vert
_{\infty }.
\end{multline}
\end{theorem}

\begin{proof}
Let $x,t\in \left[ a,b\right] $ with $t\neq x.$ Applying Cauchy's mean value
theorem, there exists a $\eta $ between $t$ and $x$ such that 
\begin{equation*}
\left( f\left( x\right) -f\left( t\right) \right) =\frac{f^{\prime }\left(
\eta \right) }{g^{\prime }\left( \eta \right) }\left( g\left( x\right)
-g\left( t\right) \right)
\end{equation*}
from where we get 
\begin{equation}
\left| f\left( x\right) -f\left( t\right) \right| =\left| \frac{f^{\prime
}\left( \eta \right) }{g^{\prime }\left( \eta \right) }\right| \left|
g\left( x\right) -g\left( t\right) \right| \leq \left\| \frac{f^{\prime }}{%
g^{\prime }}\right\| _{\infty }\left| g\left( x\right) -g\left( t\right)
\right| ,  \label{2.3}
\end{equation}
for any $t,x\in \left[ a,b\right] .$

Using the properties of the integral, we deduce by (\ref{2.3}), that 
\begin{align}
\left| f\left( x\right) -\frac{1}{b-a}\int_{a}^{b}f\left( t\right) dt\right|
& \leq \frac{1}{b-a}\int_{a}^{b}\left| f\left( x\right) -f\left( t\right)
\right| dt  \label{2.4} \\
& \leq \left\| \frac{f^{\prime }}{g^{\prime }}\right\| _{\infty }\frac{1}{b-a%
}\int_{a}^{b}\left| g\left( x\right) -g\left( t\right) \right| dt.  \notag
\end{align}
Since $g^{\prime }\left( t\right) \neq 0$ on $\left( a,b\right) ,$ it
follows that either $g^{\prime }\left( t\right) >0$ or $g^{\prime }\left(
t\right) <0$ for any $t\in \left( a,b\right) .$

If $g^{\prime }\left( t\right) >0$ for all $t\in \left( a,b\right) ,$ then $%
g $ is strictly monotonic increasing on $\left( a,b\right) $ and 
\begin{align*}
\int_{a}^{b}\left| g\left( x\right) -g\left( t\right) \right| dt&
=\int_{a}^{x}\left( g\left( x\right) -g\left( t\right) \right)
dt+\int_{x}^{b}\left( g\left( t\right) -g\left( x\right) \right) dt \\
& =2\left( x-\frac{a+b}{2}\right) g\left( x\right) +\int_{x}^{b}g\left(
t\right) dt-\int_{a}^{x}g\left( t\right) dt.
\end{align*}
If $g^{\prime }\left( t\right) <0$ for all $t\in \left( a,b\right) ,$ then 
\begin{equation*}
\int_{a}^{b}\left| g\left( x\right) -g\left( t\right) \right| dt=-\left[
2\left( x-\frac{a+b}{2}\right) g\left( x\right) +\int_{x}^{b}g\left(
t\right) dt-\int_{a}^{x}g\left( t\right) dt\right]
\end{equation*}
and the inequality (\ref{2.2}) is proved.
\end{proof}

The following midpoint inequality is a natural consequence of the above
result.

\begin{corollary}
\label{c2.2}With the above assumptions for $f$ and $g,$ one has the
inequality 
\begin{multline}
\left| f\left( \frac{a+b}{2}\right) -\frac{1}{b-a}\int_{a}^{b}f\left(
t\right) dt\right|  \label{2.5} \\
\leq \frac{1}{b-a}\left| \int_{\frac{a+b}{2}}^{b}g\left( t\right)
dt-\int_{a}^{\frac{a+b}{2}}g\left( t\right) dt\right| \left\| \frac{%
f^{\prime }}{g^{\prime }}\right\| _{\infty }.
\end{multline}
\end{corollary}

\begin{remark}
\begin{enumerate}
\item If in the above theorem, we choose $g\left( t\right) =t,$ then from (%
\ref{2.2}) we recapture Ostrowski's inequality (\ref{1.1}).

\item If in Theorem \ref{t2.1} we choose $g\left( t\right) =t^{p}$, $p\in 
\mathbb{R}\backslash \left\{ 0\right\} ,$ or $g\left( t\right) =\ln t$ with $%
t\in \left( a,b\right) \subset \left( 0,\infty \right) ,$ then we obtain
Theorem \ref{t2} and Theorem \ref{t3} respectively.
\end{enumerate}
\end{remark}

One may obtain many inequalities from Theorem \ref{t2.1} on choosing
different instances of functions $g.$

\begin{proposition}
\label{p2.3}Let $f:\left[ a,b\right] \subset \mathbb{R}\rightarrow \mathbb{R}
$ be continuous on $\left[ a,b\right] $ and differentiable on $\left(
a,b\right) .$ If there exists a constant $\Gamma <\infty $ such that 
\begin{equation}
\left\vert f^{\prime }\left( t\right) \right\vert \leq \Gamma e^{-t}\text{ 
\hspace{0.05in}for any \hspace{0.05in}}t\in \left( a,b\right) ,  \label{2.6}
\end{equation}%
then one has the inequality: 
\begin{multline}
\left\vert f\left( \frac{a+b}{2}\right) -\frac{1}{b-a}\int_{a}^{b}f\left(
t\right) dt\right\vert  \label{2.7} \\
\leq \Gamma \left[ 2\left( \frac{x-A\left( a,b\right) }{b-a}\right) e^{x}+%
\frac{\left( b-x\right) E\left( x,b\right) -\left( x-a\right) E\left(
a,x\right) }{b-a}\right]
\end{multline}%
for any $x\in \left( a,b\right) ,$ where $A=A\left( a,b\right) =\frac{a+b}{2}
$ and $E$ is the exponential men, i.e., 
\begin{equation*}
E\left( x,y\right) :=\left\{ 
\begin{array}{ll}
\dfrac{e^{x}-e^{y}}{x-y} & \text{if \hspace{0.05in}}x\neq y \\ 
&  \\ 
e^{y} & \text{if \hspace{0.05in}}x=y%
\end{array}%
\right. ,\;x,y\in \mathbb{R}\text{.}
\end{equation*}%
In particular, we have 
\begin{equation}
\left\vert f\left( A\right) -\frac{1}{b-a}\int_{a}^{b}f\left( t\right)
dt\right\vert \leq \frac{1}{2}\left[ E\left( A,b\right) -E\left( a,A\right) %
\right] \Gamma .  \label{2.8}
\end{equation}
\end{proposition}

The proof is obvious by Theorem \ref{t2.1} on choosing $g\left( t\right)
=e^{t}$ and we omit the details.

Another example is considered in the following proposition.

\begin{proposition}
\label{p2.4}Let $f:\left[ a,b\right] \subset \left( 0,\frac{\pi }{2}\right)
\rightarrow \mathbb{R}$ be continuous on $\left[ a,b\right] $ and
differentiable on $\left( a,b\right) .$

\begin{enumerate}
\item[$\left( i\right) $] If there exists a constant $\Gamma _{1}<\infty $ 
\hspace{0.05in}such that 
\begin{equation}
\left| f^{\prime }\left( t\right) \right| \leq \Gamma _{1}\cos t,\;t\in
\left( a,b\right) ,  \label{2.9}
\end{equation}
then one has the inequality 
\begin{multline}
\left| f\left( x\right) -\frac{1}{b-a}\int_{a}^{b}f\left( t\right) dt\right|
\label{2.10} \\
\leq \Gamma _{1}\left[ 2\left( \frac{x-A\left( a,b\right) }{b-a}\right) \sin
x+\frac{\left( x-a\right) C\left( a,x\right) -\left( b-x\right) C\left(
x,b\right) }{b-a}\right]
\end{multline}
for any $x\in \left( a,b\right) ,$ where $C$ is the cos-mean value, i.e., 
\begin{equation*}
C\left( x,y\right) :=\left\{ 
\begin{array}{ll}
\dfrac{\cos x-\cos y}{x-y} & \text{if \hspace{0.05in}}x\neq y \\ 
&  \\ 
-\sin y & \text{if \hspace{0.05in}}x=y%
\end{array}
.\right.
\end{equation*}
In particular we have 
\begin{equation}
\left| f\left( A\right) -\frac{1}{b-a}\int_{a}^{b}f\left( t\right) dt\right|
\leq \frac{1}{2}\left[ C\left( a,A\right) -C\left( A,b\right) \right] \Gamma
_{1}.  \label{2.11}
\end{equation}

\item[$\left( ii\right) $] If there exists a constant $\Gamma _{2}<\infty $
such that 
\begin{equation}
\left| f^{\prime }\left( t\right) \right| \leq \Gamma _{1}\sin t,\;t\in
\left( a,b\right) ,  \label{2.12}
\end{equation}
then one has the inequality 
\begin{multline}
\left| f\left( x\right) -\frac{1}{b-a}\int_{a}^{b}f\left( t\right) dt\right|
\label{2.13} \\
\leq \Gamma _{2}\left[ 2\left( \frac{x-A\left( a,b\right) }{b-a}\right) \cos
x+\frac{\left( b-x\right) S\left( x,b\right) -\left( x-a\right) S\left(
a,x\right) }{b-a}\right] ,
\end{multline}
for any $x\in \left( a,b\right) ,$ where $S$ is the $\sin -$mean value,
i.e., 
\begin{equation*}
S\left( x,y\right) :=\left\{ 
\begin{array}{ll}
\dfrac{\sin x-\sin y}{x-y} & \text{if \hspace{0.05in}}x\neq y \\ 
&  \\ 
\cos y & \text{if \hspace{0.05in}}x=y%
\end{array}
.\right.
\end{equation*}
In particular, we have 
\begin{equation}
\left| f\left( A\right) -\frac{1}{b-a}\int_{a}^{b}f\left( t\right) dt\right|
\leq \frac{1}{2}\left[ S\left( A,b\right) -S\left( a,A\right) \right] \Gamma
_{2}.  \label{2.14}
\end{equation}
\end{enumerate}
\end{proposition}

The following result also holds.

\begin{theorem}
\label{t2.5}Let $f,g:\left[ a,b\right] \rightarrow \mathbb{R}$ be continuous
on $\left[ a,b\right] $ and differentiable on $\left( a,b\right) \backslash
\left\{ x\right\} ,$ $x\in \left( a,b\right) .$ If $g^{\prime }\left(
t\right) \neq 0$ for $t\in \left( a,x\right) \cup \left( x,b\right) ,$ then
we have the inequality 
\begin{multline}
\left\vert f\left( x\right) -\frac{1}{b-a}\int_{a}^{b}f\left( t\right)
dt\right\vert  \label{2.15} \\
\leq \frac{1}{b-a}\left\vert g\left( x\right) \left( x-a\right)
-\int_{a}^{x}g\left( t\right) dt\right\vert \cdot \left\Vert \frac{f^{\prime
}}{g^{\prime }}\right\Vert _{\left( a,x\right) ,\infty } \\
+\frac{1}{b-a}\left\vert g\left( x\right) \left( b-x\right)
-\int_{x}^{b}g\left( t\right) dt\right\vert \cdot \left\Vert \frac{f^{\prime
}}{g^{\prime }}\right\Vert _{\left( x,b\right) ,\infty }.
\end{multline}
\end{theorem}

\begin{proof}
We obviously have: 
\begin{align}
& \left\vert f\left( x\right) -\frac{1}{b-a}\int_{a}^{b}f\left( t\right)
dt\right\vert  \label{2.16} \\
& =\left\vert \frac{1}{b-a}\int_{a}^{b}\left( f\left( x\right) -f\left(
t\right) \right) dt\right\vert  \notag \\
& \leq \frac{1}{b-a}\int_{a}^{b}\left\vert f\left( x\right) -f\left(
t\right) \right\vert dt  \notag \\
& =\frac{1}{b-a}\left[ \int_{a}^{x}\left\vert f\left( x\right) -f\left(
t\right) \right\vert dt+\int_{x}^{b}\left\vert f\left( x\right) -f\left(
t\right) \right\vert dt\right] .  \notag
\end{align}%
Applying Cauchy's mean value theorem on the interval $\left( a,x\right) ,$
we deduce (see the proof of Theorem \ref{t2.1}) that 
\begin{equation}
\left\vert f\left( x\right) -f\left( t\right) \right\vert \leq \left\Vert 
\frac{f^{\prime }}{g^{\prime }}\right\Vert _{\left( a,x\right) ,\infty
}\left\vert g\left( x\right) -g\left( t\right) \right\vert  \label{2.17}
\end{equation}%
for any $t\in \left( a,x\right) ,$ and, similarly 
\begin{equation}
\left\vert f\left( x\right) -f\left( t\right) \right\vert \leq \left\Vert 
\frac{f^{\prime }}{g^{\prime }}\right\Vert _{\left( x,b\right) ,\infty
}\left\vert g\left( x\right) -g\left( t\right) \right\vert  \label{2.18}
\end{equation}%
for any $t\in \left( x,b\right) .$

Consequently 
\begin{equation*}
\int_{a}^{x}\left| f\left( x\right) -f\left( t\right) \right| dt\leq \left\| 
\frac{f^{\prime }}{g^{\prime }}\right\| _{\left( a,x\right) ,\infty
}\int_{a}^{x}\left| g\left( x\right) -g\left( t\right) \right| dt
\end{equation*}
and 
\begin{equation*}
\int_{x}^{b}\left| f\left( x\right) -f\left( t\right) \right| dt\leq \left\| 
\frac{f^{\prime }}{g^{\prime }}\right\| _{\left( x,b\right) ,\infty
}\int_{x}^{b}\left| g\left( x\right) -g\left( t\right) \right| dt.
\end{equation*}
Since $g^{\prime }$ has a constant sign in either $\left( a,x\right) $ or $%
\left( x,b\right) ,$ it follows that $g$ is strictly increasing or strictly
decreasing in $\left( a,x\right) $ and $\left( x,b\right) .$

Thus 
\begin{eqnarray*}
\int_{a}^{x}\left| g\left( x\right) -g\left( t\right) \right| dt &=&\left\{ 
\begin{array}{ll}
g\left( x\right) \left( x-a\right) -\int_{a}^{x}g\left( t\right) dt & \text{%
if }g\text{\hspace{0.05in}is increasing on }\left[ a,x\right] \\ 
&  \\ 
\int_{a}^{x}g\left( t\right) dt-g\left( x\right) \left( x-a\right) & \text{%
if }g\text{\hspace{0.05in}is decreasing}%
\end{array}
\right. \\
&=&\left| g\left( x\right) \left( x-a\right) -\int_{a}^{x}g\left( t\right)
dt\right|
\end{eqnarray*}
and, in a similar way 
\begin{equation*}
\int_{x}^{b}\left| g\left( x\right) -g\left( t\right) \right| dt=\left|
g\left( x\right) \left( b-x\right) -\int_{x}^{b}g\left( t\right) dt\right| .
\end{equation*}
Consequently, by the use of (\ref{2.16}), we deduce the desired inequality (%
\ref{2.15}).
\end{proof}

The following particular case may be of interest.

\begin{corollary}
\label{c2.6}Let $f,g:\left[ a,b\right] \rightarrow \mathbb{R}$ be continuous
on $\left[ a,b\right] $ and differentiable on $\left( a,b\right) \backslash
\left\{ \frac{a+b}{2}\right\} .$ If $g^{\prime }\left( t\right) \neq 0$ on $%
\left( a,\frac{a+b}{2}\right) \cup \left( \frac{a+b}{2},b\right) ,$ then we
have the inequality 
\begin{multline}
\left\vert f\left( \frac{a+b}{2}\right) -\frac{1}{b-a}\int_{a}^{b}f\left(
t\right) dt\right\vert  \label{2.19} \\
\leq \frac{1}{2}\left\{ \left\vert g\left( \frac{a+b}{2}\right) -\frac{2}{b-a%
}\int_{a}^{\frac{a+b}{2}}g\left( t\right) dt\right\vert \cdot \left\Vert 
\frac{f^{\prime }}{g^{\prime }}\right\Vert _{\left( a,\frac{a+b}{2}\right)
,\infty }\right. \\
+\left. \left\vert g\left( \frac{a+b}{2}\right) -\frac{2}{b-a}\int_{\frac{a+b%
}{2}}^{b}g\left( t\right) dt\right\vert \cdot \left\Vert \frac{f^{\prime }}{%
g^{\prime }}\right\Vert _{\left( \frac{a+b}{2},b\right) ,\infty }\right\} .
\end{multline}
\end{corollary}

The following result also holds.

\begin{proposition}
\label{p2.7}Let $f:\left[ a,b\right] \rightarrow \mathbb{R}$ be continuous
on $\left[ a,b\right] $ and differentiable on $\left( a,b\right) \backslash
\left\{ x\right\} ,$ $x\in \left( a,b\right) .$ Assume that, for $p>0,$ we
have 
\begin{equation}
\left| f^{\prime }\left( t\right) \right| \leq \left\{ 
\begin{array}{ll}
M_{1,p}\left( x\right) \left( x-t\right) ^{1-p} & \text{for any \hspace{%
0.05in}}t\in \left( a,x\right) , \\ 
&  \\ 
M_{2,p}\left( x\right) \left( t-x\right) ^{1-p} & \text{for any \hspace{%
0.05in}}t\in \left( x,b\right) .%
\end{array}
\right.  \label{2.20}
\end{equation}
Then we have the inequality 
\begin{multline}
\left| f\left( x\right) -\frac{1}{b-a}\int_{a}^{b}f\left( t\right) dt\right|
\label{2.21} \\
\leq \frac{1}{p\left( p+1\right) }\left( b-a\right) \left[ M_{1,p}\left(
x\right) \left( x-a\right) ^{p+1}+M_{2,p}\left( x\right) \left( b-x\right)
^{p+1}\right] .
\end{multline}
\end{proposition}

The proof follows by Theorem \ref{t2.5} applied for $g\left( x\right)
=\left| x-t\right| ^{p},$ $p>0.$ We omit the details.

\begin{remark}
\label{r2.8}If $f$ is as in Proposition \ref{p2.7} and 
\begin{equation}
\left| f^{\prime }\left( t\right) \right| \leq \left\{ 
\begin{array}{ll}
M_{1}\left( \frac{a+b}{2}\right) \left( \frac{a+b}{2}-t\right) ^{1-p} & 
\text{for any \hspace{0.05in}}t\in \left( a,\frac{a+b}{2}\right) , \\ 
&  \\ 
M_{2}\left( \frac{a+b}{2}\right) \left( t-\frac{a+b}{2}\right) ^{1-p} & 
\text{for any \hspace{0.05in}}t\in \left( \frac{a+b}{2},b\right) ,%
\end{array}
\right.  \label{2.22}
\end{equation}
then, by (\ref{2.21}), we get 
\begin{multline}
\left| f\left( \frac{a+b}{2}\right) -\frac{1}{b-a}\int_{a}^{b}f\left(
t\right) dt\right|  \label{2.23} \\
\leq \frac{\left( b-a\right) ^{p+1}}{2^{p+1}p\left( p+1\right) }\left[
M_{1}\left( \frac{a+b}{2}\right) +M_{2}\left( \frac{a+b}{2}\right) \right] .
\end{multline}
\end{remark}

\begin{remark}
\label{r2.9}If $f$ is as in Proposition \ref{p2.7} and 
\begin{equation*}
\left| f^{\prime }\left( t\right) \right| \leq M_{p}\left( x\right) \left|
x-t\right| ^{1-p}\text{\hspace{0.05in}}t\in \left( a,b\right) ,
\end{equation*}
then, by (\ref{2.21}), we get 
\begin{multline}
\left| f\left( x\right) -\frac{1}{b-a}\int_{a}^{b}f\left( t\right) dt\right|
\label{2.24} \\
\leq \frac{1}{p\left( p+1\right) \left( b-a\right) }\left[ \left( x-a\right)
^{p+1}+\left( b-x\right) ^{p+1}\right] M_{p}\left( x\right) ,
\end{multline}
which is the result obtained in (\ref{1.4}).
\end{remark}

\section{Some Inequalities of Midpoint Type}

\begin{enumerate}
\item Let $0<a<b.$ Consider the function $g:\left[ a,b\right] \rightarrow 
\mathbb{R}$, $g\left( t\right) =t^{p},$ $t\in \mathbb{R}\backslash \left\{
0,-1\right\} .$ Then $g^{\prime }\left( t\right) =pt^{p-1},$ $g\left( \frac{%
a+b}{2}\right) =A^{p}\left( a,b\right) ,$%
\begin{align*}
\frac{2}{b-a}\int_{a}^{\frac{a+b}{2}}g\left( t\right) dt& =L_{p}^{p}\left(
a,A\left( a,b\right) \right) , \\
\frac{2}{b-a}\int_{\frac{a+b}{2}}^{b}g\left( t\right) dt& =L_{p}^{p}\left(
A\left( a,b\right) ,b\right) ,
\end{align*}
and by Corollary \ref{c2.6}, we may state the following proposition.

\begin{proposition}
\label{p3.1}Let $f:\left[ a,b\right] \subset \left( 0,\infty \right)
\rightarrow \mathbb{R}$ be continuous on $\left[ a,b\right] $ and
differentiable on $\left( a,b\right) \backslash \left\{ \frac{a+b}{2}%
\right\} .$ If 
\begin{equation}
\left| f^{\prime }\left( t\right) \right| \leq \left\{ 
\begin{array}{ll}
M_{1}\left( \frac{a+b}{2}\right) t^{p}, & t\in \left( a,\frac{a+b}{2}\right)
, \\ 
&  \\ 
M_{2}\left( \frac{a+b}{2}\right) t^{p}, & t\in \left( \frac{a+b}{2},b\right)
,%
\end{array}
\right.  \label{3.0}
\end{equation}
then we have the inequality 
\begin{multline}
\left| f\left( \frac{a+b}{2}\right) -\frac{1}{b-a}\int_{a}^{b}f\left(
t\right) dt\right|  \label{3.1} \\
\leq \frac{1}{2p}\left\{ M_{1}\left( \frac{a+b}{2}\right) \left| A^{p}\left(
a,b\right) -L_{p}^{p}\left( a,A\left( a,b\right) \right) \right| \right. \\
+\left. M_{2}\left( \frac{a+b}{2}\right) \left| L_{p}^{p}\left( A\left(
a,b\right) ,b\right) -A^{p}\left( a,b\right) \right| \right\} .
\end{multline}
\end{proposition}

The particular case $p=1$ is of interest and so we may state the following
corollary.

\begin{corollary}
\label{c3.2}Let $f:\left[ a,b\right] \subset \left( 0,\infty \right)
\rightarrow \mathbb{R}$ be continuous on $\left[ a,b\right] $ and
differentiable on $\left( a,b\right) \backslash \left\{ \frac{a+b}{2}%
\right\} .$ If 
\begin{equation}
\left| f^{\prime }\left( t\right) \right| \leq \left\{ 
\begin{array}{ll}
N_{1}\left( \frac{a+b}{2}\right) t, & t\in \left( a,\frac{a+b}{2}\right) ,
\\ 
&  \\ 
N_{2}\left( \frac{a+b}{2}\right) t, & t\in \left( \frac{a+b}{2},b\right) ,%
\end{array}
\right.  \label{3.2}
\end{equation}
then we have the inequality: 
\begin{multline}
\left| f\left( \frac{a+b}{2}\right) -\frac{1}{b-a}\int_{a}^{b}f\left(
t\right) dt\right|  \label{3.3} \\
\leq \frac{1}{8}\left[ N_{1}\left( \frac{a+b}{2}\right) +N_{2}\left( \frac{%
a+b}{2}\right) \right] \left( b-a\right) .
\end{multline}
\end{corollary}

\item Let $0<a<b.$ Consider the function $g:\left[ a,b\right] \rightarrow 
\mathbb{R}$, $g\left( t\right) =\frac{1}{t}.$ Then $g^{\prime }\left(
t\right) =-\frac{1}{t^{2}},$ $g\left( \frac{a+b}{2}\right) =A^{-1}\left(
a,b\right) ,$%
\begin{align*}
\frac{2}{b-a}\int_{a}^{\frac{a+b}{2}}g\left( t\right) dt& =L^{-1}\left(
a,A\left( a,b\right) \right) , \\
\frac{2}{b-a}\int_{\frac{a+b}{2}}^{b}g\left( t\right) dt& =L^{-1}\left(
A\left( a,b\right) ,b\right) ,
\end{align*}
and by Corollary \ref{c2.6} we may state the following Proposition.

\begin{proposition}
\label{p3.3}Let $f:\left[ a,b\right] \subset \left( 0,\infty \right)
\rightarrow \mathbb{R}$ be continuous on $\left[ a,b\right] $ and
differentiable on $\left( a,b\right) \backslash \left\{ \frac{a+b}{2}%
\right\} .$ If 
\begin{equation}
\left| f^{\prime }\left( t\right) \right| \leq \left\{ 
\begin{array}{ll}
M_{1}\left( \frac{a+b}{2}\right) t^{-2}, & t\in \left( a,\frac{a+b}{2}%
\right) , \\ 
&  \\ 
M_{2}\left( \frac{a+b}{2}\right) t^{-2}, & t\in \left( \frac{a+b}{2}%
,b\right) ,%
\end{array}
\right.  \label{3.4}
\end{equation}
then we have the inequality: 
\begin{multline}
\left| f\left( \frac{a+b}{2}\right) -\frac{1}{b-a}\int_{a}^{b}f\left(
t\right) dt\right|  \label{3.5} \\
\leq \frac{1}{2}\left[ M_{1}\left( \frac{a+b}{2}\right) \cdot \frac{\left[
A\left( a,b\right) -L\left( a,A\left( a,b\right) \right) \right] }{L\left(
a,A\left( a,b\right) \right) A\left( a,b\right) }\right. \\
+\left. M_{2}\left( \frac{a+b}{2}\right) \cdot \frac{\left[ L\left( A\left(
a,b\right) ,b\right) -A\left( a,b\right) \right] }{L\left( A\left(
a,b\right) ,b\right) A\left( a,b\right) }\right] .
\end{multline}
\end{proposition}

\item Let $0<a<b.$ Consider the function $g:\left[ a,b\right] \rightarrow 
\mathbb{R}$, $g\left( t\right) =\ln t.$ Then $g^{\prime }\left( t\right) =%
\frac{1}{t}$, $g\left( \frac{a+b}{2}\right) =\ln A\left( a,b\right) ,$%
\begin{align*}
\frac{2}{b-a}\int_{a}^{\frac{a+b}{2}}g\left( t\right) dt& =\ln I\left(
a,A\left( a,b\right) \right) , \\
\frac{2}{b-a}\int_{\frac{a+b}{2}}^{b}g\left( t\right) dt& =\ln I\left(
A\left( a,b\right) ,b\right) ,
\end{align*}
and by Corollary \ref{c2.6} we may state the following proposition.

\begin{proposition}
\label{p3.4}Let $f:\left[ a,b\right] \subset \left( 0,\infty \right)
\rightarrow \mathbb{R}$ be continuous on $\left[ a,b\right] $ and
differentiable on $\left( a,b\right) \backslash \left\{ \frac{a+b}{2}%
\right\} .$ If 
\begin{equation}
\left| f^{\prime }\left( t\right) \right| \leq \left\{ 
\begin{array}{ll}
M_{1}\left( \frac{a+b}{2}\right) t^{-1}, & t\in \left( a,\frac{a+b}{2}%
\right) , \\ 
&  \\ 
M_{2}\left( \frac{a+b}{2}\right) t^{-1}, & t\in \left( \frac{a+b}{2}%
,b\right) ,%
\end{array}
\right.  \label{3.6}
\end{equation}
then we have the inequality: 
\begin{multline}
\left| f\left( \frac{a+b}{2}\right) -\frac{1}{b-a}\int_{a}^{b}f\left(
t\right) dt\right|  \label{3.7} \\
\leq \ln \left\{ G\left( \left[ \frac{A\left( a,b\right) }{I\left( a,A\left(
a,b\right) \right) }\right] ^{M_{1}\left( \frac{a+b}{2}\right) },\left[ 
\frac{I\left( A\left( a,b\right) ,b\right) }{A\left( a,b\right) }\right]
^{M_{2}\left( \frac{a+b}{2}\right) }\right) \right\} .
\end{multline}
\end{proposition}
\end{enumerate}

\section{The Case of Weighed Integrals}

We may state the following theorem.

\begin{theorem}
\label{t4.1}Let $f,g:\left[ a,b\right] \rightarrow \mathbb{R}$ be continuous
on $\left[ a,b\right] $ and differentiable on $\left( a,b\right) $ and $w:%
\left[ a,b\right] \rightarrow \lbrack 0,\infty )$ an integrable function
such that $\int_{a}^{b}w\left( s\right) ds>0.$ If $g^{\prime }\left(
t\right) \neq 0$ for each $t\in \left( a,b\right) $ and 
\begin{equation}
\left\Vert \frac{f^{\prime }}{g^{\prime }}\right\Vert _{\infty
}:=\sup\limits_{t\in \left( a,b\right) }\left\vert \frac{f^{\prime }\left(
t\right) }{g^{\prime }\left( t\right) }\right\vert <\infty ,  \label{4.1}
\end{equation}%
then for any $x\in \left( a,b\right) $ one has the inequality 
\begin{multline}
\left\vert f\left( x\right) -\frac{1}{\int_{a}^{b}w\left( t\right) dt}%
\int_{a}^{b}f\left( t\right) w\left( t\right) dt\right\vert  \label{4.2} \\
\leq \left\vert g\left( x\right) \cdot \frac{\int_{a}^{x}w\left( t\right)
dt-\int_{x}^{b}w\left( t\right) dt}{\int_{a}^{b}w\left( t\right) dt}+\frac{%
\int_{x}^{b}w\left( t\right) g\left( t\right) dt-\int_{a}^{x}g\left(
t\right) w\left( t\right) dt}{\int_{a}^{b}w\left( t\right) dt}\right\vert
\cdot \left\Vert \frac{f^{\prime }}{g^{\prime }}\right\Vert _{\infty }.
\end{multline}
\end{theorem}

\begin{proof}
Let $x,t\in \left[ a,b\right] $ with $t\neq x.$ Applying Cauchy's mean value
theorem, there exists a $\eta $ between $t$ and $x$ such that 
\begin{equation*}
f\left( x\right) -f\left( t\right) =\frac{f^{\prime }\left( \eta \right) }{%
g^{\prime }\left( \eta \right) }\left[ g\left( x\right) -g\left( t\right) %
\right] ,
\end{equation*}
from where we get 
\begin{equation}
\left| f\left( x\right) -f\left( t\right) \right| =\left| \frac{f^{\prime
}\left( \eta \right) }{g^{\prime }\left( \eta \right) }\right| \left|
g\left( x\right) -g\left( t\right) \right| \leq \left\| \frac{f^{\prime }}{%
g^{\prime }}\right\| _{\infty }\left| g\left( x\right) -g\left( t\right)
\right|  \label{4.3}
\end{equation}
for any $t,x\in \left[ a,b\right] .$

Using the properties of the integral, we deduce by (\ref{4.3}), that 
\begin{align}
& \left| f\left( x\right) -\frac{1}{\int_{a}^{b}w\left( s\right) ds}%
\int_{a}^{b}w\left( s\right) f\left( s\right) ds\right|  \label{4.4} \\
& \leq \frac{1}{\int_{a}^{b}w\left( s\right) ds}\int_{a}^{b}w\left( t\right)
\left| f\left( x\right) -f\left( t\right) \right| dt  \notag \\
& \leq \left\| \frac{f^{\prime }}{g^{\prime }}\right\| _{\infty }\frac{1}{%
\int_{a}^{b}w\left( s\right) ds}\int_{a}^{b}w\left( t\right) \left| g\left(
x\right) -g\left( t\right) \right| dt.  \notag
\end{align}
Since $g^{\prime }\left( t\right) \neq 0$ on $\left( a,b\right) ,$ it
follows that either $g^{\prime }\left( t\right) >0$ or $g^{\prime }\left(
t\right) <0$ for any $t\in \left( a,b\right) .$

If $g^{\prime }\left( t\right) >0$ for all $t\in \left( a,b\right) ,$ then $%
g $ is strictly monotonic increasing on $\left( a,b\right) $ and 
\begin{align*}
& \int_{a}^{b}w\left( t\right) \left| g\left( x\right) -g\left( t\right)
\right| dt \\
& =\int_{a}^{x}w\left( t\right) \left( g\left( x\right) -g\left( t\right)
\right) dt+\int_{x}^{b}w\left( t\right) \left( g\left( t\right) -g\left(
x\right) \right) dt \\
& =g\left( x\right) \int_{a}^{x}w\left( t\right) dt-\int_{a}^{x}w\left(
t\right) g\left( t\right) dt+\int_{x}^{b}w\left( t\right) g\left( t\right)
dt-g\left( x\right) \int_{x}^{b}w\left( t\right) dt \\
& =g\left( x\right) \left[ \int_{a}^{x}w\left( t\right)
dt-\int_{x}^{b}w\left( t\right) dt\right] +\int_{x}^{b}w\left( t\right)
g\left( t\right) dt-\int_{a}^{x}w\left( t\right) g\left( t\right) dt.
\end{align*}

If $g^{\prime }\left( t\right) <0$ for all $t\in \left( a,b\right) ,$ then 
\begin{multline*}
\int_{a}^{b}w\left( t\right) \left\vert g\left( x\right) -g\left( t\right)
\right\vert dt \\
=-\left[ g\left( x\right) \left[ \int_{a}^{x}w\left( t\right)
dt-\int_{x}^{b}w\left( t\right) dt\right] +\int_{x}^{b}w\left( t\right)
g\left( t\right) dt-\int_{a}^{x}w\left( t\right) g\left( t\right) dt\right] ,
\end{multline*}%
and the inequality (\ref{2.2}) is proved.
\end{proof}

\begin{corollary}
\label{c4.2}If $x_{0}\in \left[ a,b\right] $ is a point for which 
\begin{equation}
\int_{a}^{x_{0}}w\left( t\right) dt=\int_{x_{0}}^{b}w\left( t\right) dt,
\label{4.5}
\end{equation}%
and $f,g,w$ are as in Theorem \ref{t4.1}, then we have the inequality 
\begin{multline}
\left\vert f\left( x_{0}\right) -\frac{1}{\int_{a}^{b}w\left( t\right) dt}%
\int_{a}^{b}w\left( t\right) f\left( t\right) dt\right\vert   \label{4.6} \\
\leq \frac{\left\vert \int_{x_{0}}^{b}w\left( t\right) g\left( t\right)
dt-\int_{a}^{x_{0}}g\left( t\right) w\left( t\right) dt\right\vert }{%
\int_{a}^{b}w\left( t\right) dt}\cdot \left\Vert \frac{f^{\prime }}{%
g^{\prime }}\right\Vert _{\infty }.
\end{multline}
\end{corollary}

In a similar manner, we may deduce the following result as well.

\begin{theorem}
\label{t4.3}Let $f,g:\left[ a,b\right] \rightarrow \mathbb{R}$ be continuous
on $\left[ a,b\right] $ and differentiable on $\left( a,b\right) \backslash
\left\{ x\right\} ,$ $x\in \left( a,b\right) .$ If $w:\left[ a,b\right]
\rightarrow \lbrack 0,\infty )$ is integrable and $\int_{a}^{b}w\left(
s\right) ds>0$ and $g^{\prime }\left( t\right) \neq 0$ for $t\in \left(
a,x\right) \cup \left( x,b\right) ,$ then we have the inequality 
\begin{multline}
\left\vert f\left( x\right) -\frac{1}{\int_{a}^{b}w\left( t\right) dt}%
\int_{a}^{b}w\left( t\right) f\left( t\right) dt\right\vert  \label{4.7} \\
\leq \left\vert g\left( x\right) \cdot \frac{\int_{a}^{x}w\left( t\right) dt%
}{\int_{a}^{b}w\left( t\right) dt}-\frac{\int_{a}^{x}w\left( t\right)
g\left( t\right) dt}{\int_{a}^{b}w\left( t\right) dt}\right\vert \cdot
\left\Vert \frac{f^{\prime }}{g^{\prime }}\right\Vert _{\left( a,x\right)
,\infty } \\
+\left\vert g\left( x\right) \cdot \frac{\int_{x}^{b}w\left( t\right) dt}{%
\int_{a}^{b}w\left( t\right) dt}-\frac{\int_{x}^{b}w\left( t\right) g\left(
t\right) dt}{\int_{a}^{b}w\left( t\right) dt}\right\vert \cdot \left\Vert 
\frac{f^{\prime }}{g^{\prime }}\right\Vert _{\left( x,b\right) ,\infty }.
\end{multline}
\end{theorem}


\begin{thebibliography}{9}
\bibitem{AO} A. OSTROWSKI, \"{U}ber die Absolutabweichung einer
differentienbaren Funktionen von ihren Integralmittelwert, \textit{Comment.
Math. Hel, }\textbf{10 }(1938), 226-227.

\bibitem{SSD1} S.S. DRAGOMIR, Some new inequalities of Ostrowski type, 
\textit{RGMIA\ Res. Rep. Coll.}, \textbf{5}(2002), Supplement, Article 11, 
\texttt{[ON\ LINE:http://rgmia.vu.edu.au/v5(E).html]}to apper in \textit{%
Gazette, Austral. Math. Soc.}

\bibitem{GA1} G. A. ANASTASSIOU, Multidimensional Ostrowski inequalities,
revisited. \textit{Acta Math. Hungar}., \textbf{97} (2002), no. 4, 339--353.

\bibitem{GA2} G. A. ANASTASSIOU, Univariate Ostrowski inequalities,
revisited. \textit{Monatsh. Math.,} \textbf{135} (2002), no. 3, 175--189.

\bibitem{SSD2} S.S. DRAGOMIR and T.M. RASSIAS (Eds),\textit{\ Ostrowski Type
Inequalities and Applications in Numerical Integration}, Kluwer Academic
Publishers, Dordrecht/Boston/London, 2002.
\end{thebibliography}
\end{document}